\documentstyle{amsart}
\input{amssym.def}
\input{amssym.tex}

\newtheorem{theorem}{Theorem}[section]
\newtheorem{lemma}{Lemma}[section]
\newtheorem{corollary}{Corollary}[section]

\theoremstyle{definition}

\newtheorem{definition}{Definition}[section]

\newtheorem{example}{Example}[section]

\theoremstyle{remark}

\def\ord{{\rm{ord}}}
\def\Cal{\cal}
\def\rank{{\rm{rank}}}
\def\ann{{\rm{Ann}}}
\def\ass{{\rm{Ass}}}
\def\gcd{{\rm{gcd}}}
\def\adj{{\rm{adj}}}
\def\fix{{\rm{Fix}}}

\begin{document}
\bibliographystyle{plain}
\title[Finitely presented algebraic dynamical systems]
{Fitting ideals for finitely presented\\ algebraic dynamical systems}
\author{M. Einsiedler}
\address{\hskip-.05in Mathematisches Institut\\
Universit{\"a}t Wien\\
Strudlhofgasse 4\\
A-1090 Vienna\\
Austria}
\email{\hskip-.05in manfred@@nelly.mat.univie.ac.at}
\author{T. Ward}
\address{\hskip-.05in School of Mathematics\\
University of East Anglia\\
Norwich NR4 7TJ\\
U.K.}
\email{\hskip-.05in t.ward@@uea.ac.uk}
\date{July 20 1998}
\thanks{The first author gratefully acknowledges
support from a London Mathematical Society Scheme 4 grant
and the hospitality of the University of East Anglia
where this work was done.}
\begin{abstract}{We consider a class of algebraic
dynamical systems introduced by Kitchens and Schmidt.
Under a weak finiteness condition -- the Descending Chain
Condition -- the dual modules have finite presentations. Using methods
from commutative algebra we show how the dynamical properties
of the system may be deduced from the Fitting ideals of
a finite free resolution of the finitely presented module.
The entropy and expansiveness
are shown to depend only on the first Fitting
ideal (and certain multiplicity data) which gives an easy
computation: in particular, no syzygy modules need
to be computed.

For ``square'' presentations (in which the number of
generators is equal to the number of relations) all the
dynamics is visible in the first Fitting ideal and
certain multiplicity data, and we show how the dynamical
properties and periodic point behaviour may be deduced
from the determinant of the matrix of relations.
}
\end{abstract}
\maketitle

\section{Introduction}

A natural family of measure--preserving
$\Bbb Z^d$--actions are provided by
commuting automorphisms of compact abelian groups.
Such actions are amenable to analysis using methods
from commutative algebra and commutative harmonic analysis.
The resulting theory, described in the papers
\cite{kitchens-schmidt-1989-main},
\cite{lind-schmidt-ward-1990},
\cite{rudolph-schmidt-1995},
\cite{schmidt-1990},
\cite{schmidt-ward-1993-schlickewei},
and the
monograph \cite{schmidt-1995-algebraic-dynamical-systems}
associates to such a dynamical system a
module over the ring of Laurent polynomials in
$d$ variables with integer coefficients, and then relates
various dynamical properties of the action to algebraic
or geometric properties of the corresponding module.
This singles out for attention the class of systems corresponding
to Neotherian modules (the systems satisfying the
{\sl Descending Chain Condition} of Kitchens and Schmidt,
\cite{kitchens-schmidt-1989-main}) and raises the
problem of computing the set of associated primes
of such modules.

Our purpose here is to exploit standard methods from
commutative algebra to study the dynamical systems corresponding
to Noetherian modules described via a {\sl finite presentation}.
Before describing this we recall the algebraic
description of such actions in \cite{kitchens-schmidt-1989-main}.
Let $R=\Bbb Z[u_1^{\pm1},\dots ,u_d^{\pm1}]$ be the ring of Laurent
polynomials with integral coefficients in the commuting variables
$u_1,\dots ,u_d$. If $\alpha$ is a $\Bbb Z^d$--action
by automorphisms of the
compact, abelian group $X$, then the dual
(character) group $M=\hat X$ of $X$
is an $R$-module under the dual $R$--action
$$
f\cdot a=\sum_{{\bold m}\in\Bbb Z^d}c_f({\bold m})\beta_{{\bold m}}(a)
$$
for all $a\in M$ and $f=\sum_{{\bold m}\in\Bbb Z^d}c_f({\bold m})u^{\bold
m}\in R$, where $u^{{\bold n}}=u_1^{n_1}\cdots u_d^{n_d}$ for every
${\bold n}=(n_1,\dots ,n_d)\in\Bbb Z^d$, and where $\beta_{\bold
n}=\widehat{\alpha_{{\bold n}}}$ is the automorphism of $M=\hat X$ dual
to $\alpha_{{\bold n}}$. In particular,
$$
\widehat{\alpha_{{\bold n}}}(a)
=\beta_{{\bold n}}(a)=u^{{\bold n}}\cdot a
$$
for all $\bold
n\in\Bbb Z^d$ and $a\in M$. Conversely, if $M$ is an
$R$-module, and
$$
\beta_{{\bold n}}^{M}(a)=u^{{\bold n}}\cdot
a
$$
for every ${\bold n}\in\Bbb Z^d$ and $a\in M$, then we obtain
a $\Bbb Z^d$-action
$$
\alpha^{M}:{\bold n}\to \alpha_{{\bold n}}^{ 
M}=\widehat{\beta_{{\bold n}}^{M}}
$$
on the compact, abelian
group $$X^{M}=\widehat{M}$$ dual to the $\Bbb
Z^d$-action $\beta^{M}:{\bold n}\to \beta_{{\bold n}}^{M}$ on
$M$.

The dynamical system $\alpha^{M}$ on $X_{M}$
satisfies the Descending Chain Condition (any decreasing sequence
of closed $\alpha^{M}$--invariant subgroups of $X_{M}$
stabilizes) if and only if the $R$--module is
Noetherian by Theorem 11.4 in \cite{kitchens-schmidt-1989-main}.
We assume from now on that $M$ is a Noetherian module,
in which case it has a finite presentation of the form
\begin{equation}
\label{finitepresentation}
M=M_A \cong
R^k/AR^n,
\end{equation}
where $M_A$ is generated as an $R$--module by
a subset with $k$ elements and the $k\times n$ matrix
$A$ defines the various relations in $M_A$.
Since free modules are not very interesting, we assume
that the rank of $A$ is $k$. If this is not the case,
then $M_A$ has a free submodule $  L$ with the
property that
$M_A/  L$ has a finite presentation in the
form (\ref{finitepresentation}) with $\rank(A)=k$.

In accordance with the spirit of the monograph
\cite{schmidt-1995-algebraic-dynamical-systems},
we would like then to be able to describe the dynamical
properties of the $\Bbb Z^d$--action $\alpha^{M_A}$
in terms of the matrix $A$.
Roughly speaking, we are able to (describe how to) compute
all the associated primes of $M_A$ from $A$ using
Auslander--Buchsbaum theory. This is enough to describe --
in principle -- the dynamical properties of $\alpha^{M_A}$.
For the special case $k=n$, or more
generally, of {\sl principal} associated primes,
we are also able to find the multiplicities of
the various associated primes, which allows the
entropy of $\alpha^{M_A}$ to be computed.
This means in particular that the entropy of $\alpha^{M_A}$
can be computed, and the
expansiveness of $\alpha^{M_A}$ can be decided,
without computing any syzygy modules.

Methods taken from commutative algebra are standard and
may all be found for example in Eisenbud's book
\cite{eisenbud}. We are grateful to Prof. Rodney Sharp for
pointing us to the right part of \cite{eisenbud}.

By ``entropy'' we mean topological entropy, as defined in Section 13 of
Schmidt's monograph \cite{schmidt-1995-algebraic-dynamical-systems}.

\section{Language from commutative algebra}

Let $S$ be a commutative ring
(recall that $R$ is the ring of Laurent polynomials
in $d$ variables with integer coefficients).
The basic terminology for an $S$--module $M$ may be found in any
commutative algebra book. A prime ideal $P\subset S$
is {\sl associated} to $M$ if there is an element $m\in M$
with the property that
\begin{equation}
\label{associatedprime}
P=\ann_{M}(m)=\{f\in S\mid
f\cdot m=0\in M\}.
\end{equation}
The module $M$ is {\sl Noetherian} if each submodule is
finitely generated (the ring $S$ is
Noetherian if it is a Noetherian $S$--module),
and this holds for Noetherian rings if and only if $M$ has
a finite presentation (\ref{finitepresentation}).
The set $\ass(M)$ of associated primes of a Noetherian
module is finite (see Theorem 6.5 in 
\cite{matsumura-1986-commutative}).
A Noetherian module is {\sl free} if it has a
presentation (\ref{finitepresentation}) in which
the matrix $A$ comprises zeros, and is {\sl cyclic}
if it has a presentation (\ref{finitepresentation})
with $k=1$. A {\sl finite free resolution}
of a Noetherian module $M$ is an exact sequence
of $S$--modules and $S$--module homomorphisms
\begin{equation}
\label{ffr}
0\longrightarrow F_n\overset{\phi_n}{\longrightarrow}\dots
\overset{\phi_2}{\longrightarrow}
F_1\overset{\phi_1}{\longrightarrow}  F_0
\longrightarrow M\longrightarrow0
\end{equation}
in which each $  F_i$ is a free $  S$--module.

A subset $U$ in $S$ is {\sl multiplicative}
if it is closed under multiplication. Each multiplicative
subset $U\subset S$ defines a {\sl localization}
\begin{equation}
\label{localization}
S^{U}=\{\textstyle\frac{s}{u}\mid s\in  S,u\in U\},
\end{equation}
where two fractions $\frac{s}{u}$ and $\frac{s'}{u'}$ are identified if there
is an element $u''\in U$ with $u''(u's-us')=0$.
The notation is altered for one special case:
if $P$ is a prime ideal in $S$, then
write $S^{(P)}$ for
$S^{S\backslash P}$, called the
localization at the prime $P$.
For a module $M$, the same definition as
(\ref{localization}) works and defines a localized
module $M^{S}$ or $M^{P}$.
If the ideal $P=\langle\pi\rangle$ is principal,
write $M^{(\pi)}=M^{P}$.
The {\sl dimension} $\dim(S)$
of $S$ is the supremum of
the length of chains of distinct prime ideals in
$S$, and this coincides with the
supremum of $\dim(S^{P})$ over all
prime ideals $P$. The dimension of a localization
$S^{P}$ is also known as the
codimension of $P$, and coincides with the
supremum of lengths of chains of prime ideals
descending from $P$.

A ring is {\sl Noetherian} if every
ascending chain of ideals stabilizes,
is a {\sl local ring} if it has just one
maximal ideal, and is {\sl regular} if 
it is Noetherian and the localization at every prime
ideal is a regular local ring. A local ring is a regular local ring if the
maximal ideal is generated by
exactly $d$ elements where $d$ is the dimension of the
local ring. It is clear that $R$
is a regular ring, and it follows (see Chapter
19 of \cite{eisenbud}) that any Noetherian $R$--module
has a finite free resolution (\ref{ffr}).
Notice that the presentation (\ref{finitepresentation})
is itself the start of a finite free resolution of
$M_A$:
$$
\dots\longrightarrow R^n\overset{A}{\longrightarrow}
R^k\longrightarrow M_A\longrightarrow0.
$$

If $M$ is a Noetherian $R$--module with
associated primes $\ass(M)=\{P_1,\dots,P_r\}$,
then there is a {\sl prime filtration} of $M$,
\begin{equation}
\label{primefiltration}
M=M_{\ell}\supset M_{\ell-1}\supset\dots
\supset M_1\supset M_0=\{0\},
\end{equation}
in which each quotient $M_j/M_{j-1}\cong
R/Q_{j}$ for some prime $Q_j\supset
P_i$ for some $i$ (see for example
Corollary 2.2 in \cite{schmidt-1990}).
The number of times a given prime $P_i$ appears
(that is, the number of $j$ for which
$Q_j=P_i$) is the {\sl multiplicity} of
$P_i$ in the filtration (\ref{primefiltration}).
If the prime ideal in question is principal and the module
$M$ has no free submodules, then the multiplicity
with which $P_i$ appears is independent of the filtration,
and we will therefore speak of the multiplicity of $P_i$
in $M$ (see Proposition 6.10 in \cite{lind-schmidt-ward-1990}).

\section{Dynamical properties}

Let $M$ be any countable $R$--module, with associated $\Bbb Z^d$--action
$\alpha^M$ on $X_M=\widehat{M}$. The following result shows how
the dynamical properties of $\alpha^M$ may be deduced from
the associated primes $\ass(M)$ of $M$. All these results
are in \cite{schmidt-1995-algebraic-dynamical-systems}; we state
them here for completeness. A {\sl generalized cyclotomic}
polynomial is an element of $R$ of the form
$u_1^{n_1}\dots u_d^{n_d}c(u_1^{m_1}\dots u_d^{m_d})$ for
some cyclotomic polynomial $c$ and ${\bold n}, {\bold m}\in
\Bbb Z^d$. Write $V(P)$ for the set of common zeros of
the elements of $P$ in $\Bbb C^d$.

\begin{theorem}
\label{dynamicalproperties}
The dynamical system $\alpha^M$ on $X_M$:

\noindent{\rm(a)} satisfies the Descending Chain Condition
on closed $\alpha^M$--invariant subgroups if and only
if $M$ is Noetherian;

\noindent{\rm(b)} is ergodic if and only if 
$\{\left(u_1^{n_1}\dots u_d^{n_d}\right)^k-1\mid
{\bold n}\in\Bbb Z^d\}\not\subset P$ for every
$k\ge 1$ and every $P\in\ass(M)$;

\noindent{\rm(c)} is mixing if and only if
$u_1^{n_1}\dots u_d^{n_d}-1\notin P$ for each
${\bold n}\in\Bbb Z^d\backslash\{0\}$ and every
$P\in\ass(M)$;

\noindent{\rm(d)} is mixing of all orders if and only
if either $P=pR$ for a rational prime $p$, or
$P\cap\Bbb Z=\{0\}$ and $\alpha^{R/P}$ is mixing
for every $P\in\ass(M)$;

\noindent{\rm(e)} has positive entropy if and only if
there is a $P\in\ass(M)$ that is principal and
not generated by a generalized cyclotomic polynomial;

\noindent{\rm(f)} has completely positive entropy if and
only if $\alpha^{R/P}$ has positive entropy for
every $P\in\ass(M)$;

\noindent{\rm(g)} is isomorphic to a Bernoulli shift if
and only if it has completely positive entropy;

\noindent{\rm(h)} is expansive if and only if
$M$ is Noetherian and $V(P)\cap\left(\Bbb S^1\right)^d=\emptyset$
for every $P\in\ass(M)$;

\noindent{\rm(i)} has a unique maximal measure if and
only if it has finite completely positive
entropy.

\end{theorem}

\begin{pf} For (a) see Theorem 11.4 in
\cite{kitchens-schmidt-1989-main}; (b) and (c) are
in Theorem 11.2 in \cite{kitchens-schmidt-1989-main};
(d) follows from Theorem 3.1 in \cite{schmidt-ward-1993-schlickewei}
and Theorem 3.3 in \cite{schmidt-1989};
(e), (f) and (i) are in
\cite{lind-schmidt-ward-1990};
(h) is Theorem 3.9 in \cite{schmidt-1990};
(g) is Theorem 1.1 in \cite{rudolph-schmidt-1995}.
\end{pf}

\section{Principal associated primes and entropy}

In this section we use localization 
to find the entropy of $\alpha^{M}_A$.

\begin{definition} Let $A$ be a $k\times n$ matrix
of rank $k$ over $R$. The determinental ideal
of $A$, $J_A\subset R$, is the ideal
$$
J_A=\langle f_1,\dots,f_{\binom{n}{k}}\rangle
$$
generated by all the $k\times k$ subdeterminants $\{f_1,\dots,f_{
 \binom{n}{k}}\}$ of $A$.
\end{definition}

For a polynomial $f\in R$, the {\sl logarithmic
Mahler measure} of $f$ is defined to be
\begin{equation}
\label{mahlermeasure}
m(f)=\int_0^1\dots\int_0^1
\log\vert f(e^{2\pi is_1},\dots,e^{2\pi is_d})\vert
\text{d}s_1\dots\text{d}s_d.
\end{equation}
For brevity, define $m(0)$ to be $\infty$.
Recall from \cite{lind-schmidt-ward-1990} that the entropy
of the dynamical system given by the cyclic module
$R/P$, where $P$ is a prime ideal, is
given by
\begin{equation}
\label{primeentropy}
h\left(\alpha^{R/P}\right)
=
\begin{cases}
0,&\text{if $P$ is non--principal;}\\
m(f),&\text{if $P=\langle f\rangle,
f\neq0$.}
\end{cases}
\end{equation}
More generally, since $R$ is a UFD, for any ideal
$Q\subset R$ there is a well--defined
greatest common divisor, and
$$
h(\alpha^{R/Q})=h(\alpha^{R/\gcd(Q)}),
$$
which is zero if $\gcd(Q)=\langle 1\rangle$ and
equal to $m(f)$ if $\gcd(Q)=\langle f\rangle$ (see
Lemma 4.5 in \cite{einsiedler-1997}).

For Noetherian modules, 
\begin{equation}
\label{moduleentropy}
h(\alpha^{M})=
\sum_{j=1}^{\ell}h(\alpha^{R/Q_j})
\end{equation}
where the prime ideals $Q_j$ are the
primes appearing in the filtration (\ref{primefiltration}).

\begin{theorem}
\label{entropy}The entropy of $\alpha^{M}$ is
given by
\begin{equation}
\label{entropyformula}
h\left(\alpha^{M_A}\right)=
m\left(
\gcd(J_A)\right).
\end{equation}
\end{theorem}

Before proving this, we indicate some examples.

\begin{example}
\label{examples}
(a) Taking $k=n=1$ and 
$A=[f]$ with an irreducible polynomial $f$, we recover the formula (
\ref{primeentropy}) in the cyclic case with a principal prime ideal.

\vskip.1in
\noindent(b) Taking $k=1$ and $n\ge 1$ we recover the
general cyclic case.

\vskip.1in
\noindent(c) If $k=n$ then formula (\ref{entropyformula})
simply reduces to $\det(A)$, which was shown in
Section 5 of \cite{lind-schmidt-ward-1990}.

\vskip.1in
\noindent(d) Let $k$ be an algebraic number field
with ring of integers ${\Cal O}_k$, and $f$ a
Laurent polynomial in $d$ variables
with coefficients in ${\Cal O}_k$.
The $\Bbb Z^d$--dynamical system $\beta$ dual to
multiplication by $u_1,\dots,u_d$ on the
${\Cal O}_k[u_1^{\pm1},\dots,u_d^{\pm1}]$--module
${\Cal M}={\Cal O}_k[u_1^{\pm1},\dots,u_d^{\pm1}]/\langle f\rangle$
is studied in \cite{einsiedler-1997}.
Taking an integral basis for ${\Cal O}_k$
shows that ${\Cal M}$ as an $R$--module
is of the form (\ref{finitepresentation}) with
$n=k$, and by (c) we see that
$$
h(\beta)=m(\det(A))=m\left(N_{k:\Bbb Q}f\right),
$$
recovering Theorem 3.10 in \cite{einsiedler-1997}.
\end{example}

\begin{lemma}\label{det_is_element}
 Each associated prime of $M_A$ contains $J_A$.
\end{lemma}
\begin{pf} 
 Let $B$ be a $k\times k$ subdeterminant of
 $A$. Then for any ${\bold{v}}\in R^k$,
 $\det(B)\cdot{\bold{v}}=BB^{\adj}{\bold v}$. Therefore the annihilator of
 any element of $M_A\cong {R}^k/A{R}^n$ contains $\det(B)$
 and also $J_A$.
\end{pf}
 
\begin{lemma}
\label{principals}
The principal associated primes
of $M_A$ are generated by the irreducible
factors of $\gcd(J_A)$. Moreover,
the multiplicity of each principal associated prime
in $M$ is equal to its multiplicity in
$\gcd(J_A)$.\end{lemma}

\begin{pf} 
It follows from Lemma \ref{det_is_element} that each element
of the set of principal
associated primes of $M$ contains the $k\times k$
subdeterminants of $A$.
This means that the generator of a
principal associated prime divides all the subdeterminants and is therefore
a factor of $\gcd(J_A)$.

Conversely, let $\{B_i\}$ be the set of
$k\times k$ subdeterminants
of $A$, and let
\begin{equation}
\label{factorsofsubdet}
\gcd\left(B_1,\dots,B_{\binom{n}{k}}\right)=
\pi_1^{e_1}\dots\pi_r^{e_r}
\end{equation}
be a factorization into irreducibles in $R$.
By Section 2 (\ref{primefiltration}) there is a prime filtration
\begin{equation}
\label{primefiltrationtwo}
M=M_{\ell}\supset M_{\ell-1}\supset\dots
\supset M_1\supset M_0=\{0\},
\end{equation}
with $M_j/M_{j-1}\cong R/Q_j$
with $Q_j\supset P$ for some $P\in\ass(M)$.
Localize (\ref{primefiltrationtwo}) at the prime ideal $\langle\pi_1\rangle$:
the pair $M_{j-1}\supset M_j$ localizes to
the pair $M_{j-1}^{(\pi_1)}\supset M_j^{(\pi_1)}$,
with quotient
$$
{R^{(\pi_1)}}/{Q_{j}^{(\pi_1)}}=
\begin{cases}
{R}/{\langle \pi_{1}\rangle}&\text{if $Q_j=\langle\pi_1\rangle$;}\\
0.&\text{if not.}
\end{cases}
$$
So (\ref{primefiltrationtwo}) collapses to a shortened
filtration of $R^{\pi_1}$--modules and we see that
the multiplicity of $\langle\pi_i\rangle$ in $M$
coincides with the multiplicity of $\langle\pi_i\rangle$ in $M^{(\pi_i)}$
for each $i=1,\dots,r$.

We are therefore reduced to studying the local case: let $\pi$ be any
one of the $\pi_i$'s, and write
$$
M^{(\pi)}_A=\left(R^{(\pi)}\right)^k/
A\left(R^{(\pi)}\right)^n.
$$
We can change our matrix by invertible (over $R^{(\pi)}$) elementary
row operations and this gives us an isomorphic module with the same
subdeterminants.
In $R^{(\pi)}$ define $\ord(f)=\ord_{\pi}(f)$ to be the 
number of times that $\pi$ divides into $f$.
Write $f\le g$ if $\ord(f)\le\ord(g)$, and with respect to this
partial ordering find (one of) the smallest entries in $A$.
Permute rows and columns in $A$ so that $a_{11}$ is a smallest
entry. Then $a_{i1}\ge a_{11}$ for $2\le i\le k$ so the quotient
$a_{i1}/a_{11}$ is an element of ${R}^{(\pi)}$ and we can subtract
multiples of the first row from the others to get a matrix of the form
$$
A_1=
\bmatrix a_{11}&\dots&&\\
0&\\
\vdots& &*\\
0\endbmatrix
$$
Repeat with $a_{22}$ and so on to produce a matrix of the form
$$
A_*=
\bmatrix
a_{11}	&\dots\\
0	&a_{22}	&\dots\\
\vdots	&0	&	\\
0	&\dots	&0	&a_{kk}&\dots
\endbmatrix
$$
in which each $a_{jj}$
is in turn the smallest non--zero element of the
submatrix $(a_{st})_{s,t\ge j}$.
Let $\ord(a_{jj})=e_{jj}$ for each $j$.

Now let ${\bold v}=(1,0,\dots,0)^{t}$,
so that
$$
\ann({\bold v}+A_*{R}^n)=\langle a_{11}\rangle,
$$
since the other columns of the matrix have a first component which is
divisible by $a_{11}$.
The map $f\mapsto f\cdot{\bold v}\in M^{(\pi)}$
gives a filtration
$$
{\bold v}R^{(\pi)}
\supset
\pi{\bold v}R^{(\pi)}
\supset
\dots
\supset
\pi^{(e_{11}-2)}{\bold v}R^{(\pi)}
\supset
\pi^{(e_{11}-1)}{\bold v}R^{(\pi)}
\supset
0
$$
of submodules of $M^{(\pi)}_A$.
As the same argument works for the other standard basis vectors it follows
that the multiplicity of $\pi$ in $M^{(\pi)}_A$ is
$\sum e_{jj}$.
Calculating all the subdeterminants shows that the greatest
common divisor is equal to the product $\prod_ja_{jj}
=\pi^{\sum_je_{jj}}$.
So the multiplicity of $\pi$ in $\gcd(J_A)$ is equal to the
multiplicity of the associated prime $(\pi)$ in a prime filtration of $ 
M_A$.
\end{pf}

\begin{pf*}{Proof of Theorem \ref{entropy}}
Use Lemma \ref{principals} and Section 2 to find the
principal associated primes and their
multiplicites; the result follows
by (\ref{moduleentropy}).\end{pf*}

For an ideal $P$ in $R$,
recall that $V(P)=\{
{\bold z}\in\Bbb C^d\mid
f({\bold z})=0{\ }\forall{\ }f\in P\}$ denotes
the set of common zeros of $P$.
Write $V(f)$ for $V(\langle f\rangle)$.
By Theorem \ref{dynamicalproperties},
$\alpha^{M}$ is expansive if and only if
$V(P)\cap(\Bbb S^1)^d=\emptyset$
for each associated prime $P\in\ass(M)$.

\begin{theorem}
\label{expansiveness}
Let $M_A$ be a finitely presented module
with $A$ of rank $k$. Then $\alpha^{M_A}$ is
expansive if and only if
\begin{equation}
\label{expansivesubdets}
(\Bbb S^1)^d\cap
\left(
\bigcap_{j=1,\dots,\binom{n}{k}}
V(\det(B_j))\right)=\emptyset,
\end{equation}
where $\{B_j\}$ is the set of $k\times k$
subdeterminants of $A$.
\end{theorem}

\begin{pf} Assume first that
$$
{\bold z}\in(\Bbb S^1)^d\cap
\left(
\bigcap_{j=1,\dots,\binom{n}{k}}
V(\det(B_j))\right).
$$
Assume that for every associated prime $P$
of $M_A$ there is a polynomial
$f_{P}\in P$
for which $f_{P}({\bold z})\neq 0$.
Then let $f=\prod_{P\in\ass(M)}f_{P}$
(so $f({\bold z})\neq 0$).
From a prime filtration of $M_A$ it is clear
that for some power $m$,
\begin{equation}
\label{vanishing}
f^mM=0.
\end{equation}
On the other hand, since ${\bold z}$ was
chosen to lie in the set of common zeros of all the
subdeterminants, in the ring
$$
M({\bold z})=\frac{\Bbb Z[{\bold z}^{\pm1}]^k}
{A({\bold z})\Bbb Z[{\bold z}^{\pm1}]^n}
$$
we have that all $k\times k$ subdeterminants of
$A({\bold z})$ vanish, so $\rank(A({\bold z}))<k$,
and in particular $M({\bold z})\neq 0$,
contradicting
(\ref{vanishing}).
It follows that if 
the intersection in (\ref{expansivesubdets}) contains a
point ${\bold z}$ then this point must lie in
$V(P)$ for some associated prime $P$,
showing that $\alpha^{M_A}$ is not expansive
by Theorem \ref{dynamicalproperties}.

Conversely, if $\alpha^{M_A}$ is not
expansive, then there is an associated prime
$P\in\ass(M_A)$ with
$V(P)\cap(\Bbb S^1)^d\ni{\bold z}$ say.
However the associated prime $P$ must
contain all the subdeterminants by Lemma \ref{det_is_element}.
so ${\bold z}\in\bigcap_{j=1,\dots,\binom{n}{k}}
V(\det(B_j)).$
\end{pf}

\section{The square case}

As remarked in Theorem \ref{dynamicalproperties},
various dynamical properties
of systems of the form $\alpha^{M}$ are
governed by properties of the set $\ass(M)$ of associated
primes of $M$. In this section we show that the associated
primes of a finite presentation with $k=n$ (the ``square case'')
are all visible in the determinant of the matrix of relations,
so the dynamics are as easy to deduce as in the case of
a cyclic module with a single principal associated prime.
We also calculate the periodic points because
{\sl a priori} one needs more information than
the associated primes to calculate this (see
Section 7 of \cite{lind-schmidt-ward-1990}).

\begin{lemma}
\label{squarecase}
If the finitely--presented module $M_A$ has
$k=n$ and $A$ has maximal rank, then the associated
prime ideals of $M$ are all given by irreducible factors
of $\det(A)$.
\end{lemma}

\begin{pf} Let $\det(A)=\pi_1^{e_1}\dots\pi_r^{e_r}$
be the decomposition into irreducibles.
By linear algebra over the quotient field of $R$ we know that $w\in AR^n$ if
and only if $\frac{1}{\det A}A^{\adj}w\in R^n$.
If an element ${\bold v}+AR^n$
has $\ann({\bold v}+AR^n)=P$
for some $P\in\ass(M_A)$, then
$$
P=\{f\in R\mid
\frac{f}{\det(A)}A^{\adj}{\bold v}\in R^k\}.
$$
Now in $\frac{1}{\det(A)}A^{\adj}v$ after all possible
cancellations there must be some $\pi_i$ in the
denominator (since ${\bold v}\notin AR^k$).
Let this denominator be $g$ say; then
$g$ must divide $f$ for all $f\in P$, so
$\ann({\bold v}+AR^n)=\langle g\rangle$.
As $P$ is prime the element $g$ must be irreducible.
It follows that all the associated primes of $M_A$
are principal and arise as factors of the determinant of $A$.

It is easy to see that the argument above also proves
that each irreducible
factor of $\det A$ gives an associated prime
(or use Lemma \ref{principals}) for the reverse inclusion.
\end{pf}

\begin{corollary}
The dynamical system $\alpha^{M_A}$ for a square matrix
$A$ is ergodic, mixing, mixing on a shape $F$, mixing of all orders,
$K$, if and only if the corresponding cyclic system
$\alpha^{R/\langle\det(A)\rangle}$ has the same property.
\end{corollary}

We are also able to compute directly the periodic points in
such systems. A {\sl period} for a $\Bbb Z^d$--action
$\alpha$ on $X$ is a lattice $\Lambda\subset\Bbb Z^d$ of full rank; the size
of the period is the (finite) index $\vert\Bbb Z^d/\Lambda\vert$.
The set of points of period $\Lambda$ is
$$\fix_{\Lambda}(\alpha)=\{x\in X\mid
\alpha_{\bold n}x=x{\ }\forall{\ }{\bold n}\in\Lambda\}.$$
Since the (multiplicative) dual group of $\Bbb Z^d$ is
$(\Bbb S^1)^d$, the annihilator $\Lambda^{\perp}$ of
$\Lambda$ is a subgroup of $(\Bbb S^1)^d$ with
cardinality $\vert\Bbb Z^d/\Lambda\vert$.

\begin{lemma}
\label{periodicpoints}
If $A$ is a square matrix of maximal rank, then
\[
{\rm{Fix}}_{\Lambda}\left(\alpha^{M_A}\right)
=
\begin{cases}
\infty &\text{if ${\ }\prod_{{\bold z}\in\Lambda^{\perp}}
		\vert\det(A)({\bold z})\vert =0;$}\\
\prod_{{\bold z}\in\Lambda^{\perp}}
\vert
\det(A)({\bold z})
\vert	&\text{if not.}
\end{cases}
\]
\end{lemma}

\begin{pf} For brevity we prove this for
square periods $\Lambda_n=n{\Bbb Z}^d$;
the general case is similar but notationally
unpleasant. We follow the method used in
\cite{lind-schmidt-ward-1990}, Section 7, exactly.

An element ${\bold x}\in X=\widehat{{R}^k/A{R}^k}$ is periodic
with respect
to $\Lambda_n$ if it annihilates $J(\Lambda_n){R}^k$ where $J(\Lambda_n)=
\langle u_1^n-1,\ldots,u_d^n-1\rangle$. So the periodic points are exactly
the elements in the dual group of 
\begin{equation}\label{per_quo}
{R}^k/(A{R}^k+J(\Lambda_n)^k).
\end{equation}
Therefore the number of periodic points is equal to the number
of elements in (\ref{per_quo}) whenever this quantity is
finite or is infinite if not. As ${R}/J(\Lambda_n)$ is
isomorphic to
${\Bbb Z}^{n^d}$ we see that the module (\ref{per_quo}) is isomorphic to
$({\Bbb Z}^F)^k/B({\Bbb Z}^F)^k$ where $F=\{1,\dots,n\}^d$
and $B$ is obtained from
$A$ by interpreting the variable $u_{\ell}$
as the shift of the $\ell$-th coordinate
in $F$. The number of periodic points in $X$ is now given by the
determinant of $B$ (or is infinite if $\det B=0$).
We calculate this quantity
using a suitable basis of
$({\Bbb C}^F)^k$. The elements in this
vector space have the form $(w_{\binom {\bold
 e} i})_{\binom {{\bold e}\in F} {i\in[1,k]}}$; use the basis 
$$
v_{\bold f}^j=(\delta_{ij}\omega^{f_1e_1+\cdots+f_de_d})_{\binom {\bold e} i}.
$$
where $\omega$ is a primitive $n$-th root.
In this basis the matrix $B$
becomes
$$
C=(a_{ij}(\omega^{e_1},\ldots,\omega^{e_d})\delta_{\bold e \bold f})_{\binom
 {\bold e} i \binom {\bold f} j}
$$
because the shift of the $\ell$-th coordinates in $F$ has $v_{\bold f}^j$ as
eigenvector with eigenvalue $\omega^{f_{\ell}}$.
The determinants are given by
$$
\det(B)=\det(C)=\prod_{{\bold e}\in F}\det(A)(\omega^{e_1},\ldots,
\omega^{e_d}),
$$
Because the matrix $C$ can be viewed as being of the form 
$$
\bmatrix D_{1}	&0	&0&\dots&0\\
	 0	&D_{2}	&0&\dots&0\\
	 \vdots &	& &     &\vdots \\
	 0	&&\dots   &0&D_{|F|}\\
\endbmatrix
$$
where the submatrices $D_j$ are obtained from $A$ by evaluation at $(
\omega^{e_1},\ldots,\omega^{e_d})$ for some $(e_1,\ldots,e_d)\in F$.
The determinant of such a matrix is the product of
the determinants of the submatrizes.
\end{pf}

For a lattice $\Lambda$, let $g(\Lambda)=
\min_{{\bold n}\in\Lambda\backslash
\{0\}}\{\Vert{\bold n}-{\bold 0}\Vert\}$.
The characterization of expansiveness and
Lemma \ref{periodicpoints} gives a very simple
proof of the general result that the growth rate of
periodic points coincides with the entropy for
expansive algebraic ${\Bbb Z}^d$--actions (see
Section 7 of \cite{lind-schmidt-ward-1990})
for finitely presented systems with $k=n$.

\begin{corollary}
\label{growthrate}
If $A$ is a square matrix of maximal rank, and
$V(\det(A))\cap(\Bbb S^1)^d=\emptyset$, then the growth
rate of periodic points is equal to the entropy:
$$
\lim_{g(\Lambda)\to\infty}
\frac{1}{\vert\Bbb Z^d/\Lambda\vert}
\log{\rm{Fix}}(\alpha^{M_A})=
h(\alpha^{M_A}).
$$
\end{corollary}

\section{The general case}

In this section we simply describe the appropriate 
results from commutative algebra and indicate by examples
how they may be used to compute associated primes in the
general case.

Fix the finite presentation (\ref{finitepresentation}) of
a Noetherian $R$--module $M$.
For each $R$--module map $\phi:R^a\to R^b$
define $J(\phi)$ to be the ideal generated by
the $\rank(\phi)\times\rank(\phi)$ subdeterminants
of a matrix for $\phi$. For maps $\phi$ appearing
in a finite free resolution, these ideals are
the {\sl Fitting ideals} of the module.
By convention $0\times 0$ determinants give the
trivial ideal $\langle 1\rangle$.

\begin{theorem}
\label{fitting}
Let 
\begin{equation}
0\longrightarrow F_n\overset{\phi_n}{\longrightarrow}\dots
\overset{\phi_2}{\longrightarrow}F_1\overset{\phi_1}{\longrightarrow}  F_0
\longrightarrow M_A\longrightarrow0
\end{equation}
be a finite free resolution of the $R$--module $M_A$.
Let $P$ be a prime ideal of $R$
with $\dim(R^{P})=\ell$.
Then $P\in\ass(M_A)$ if and only
if $P\supset J(\phi_{\ell}).$
\end{theorem}

\begin{pf} This is proved in Corollary 20.14 of
\cite{eisenbud} with the condition $\dim(R^{P})=\ell$
replaced by $\text{depth}(P\subset R)=\ell$
(see Chapter 18 of \cite{eisenbud} for this notion).
By Theorem 18.7 of \cite{eisenbud} we have that since
$R$ is regular (and hence Cohen--Macaulay by 
Section 18.5 of \cite{eisenbud}), $\text{depth}(P)=
\text{height}(P):=\dim(R^{P})$,
so the result follows.
\end{pf}

Notice that the first Fitting ideal $J(\phi_1)$
is exactly the ideal $J_A$ used above.

We now describe several examples to illustrate the kind
of calculations involved and some of the phenomena that may
arise.

\begin{example}
\label{finalexamples}
(a) Let $P=\langle f\rangle$
be a prime ideal. Then a finite free resolution of
$R/P$ is given by
$$
0\longrightarrow R\overset{[f]}
{\longrightarrow}R\longrightarrow R/P\longrightarrow0.
$$
By Theorem \ref{fitting}, we see that the associated primes
of $R/P$ comprise exactly $\{P\}$.

\vskip.1in
\noindent(b) Let $f$ be irreducible, and let $M=R
/\langle 2f\rangle.$ Then 
$$
0\longrightarrow R\overset{[2f]}
{\longrightarrow}R\longrightarrow M\longrightarrow0.
$$
is a free resolution of $M$.
If $\dim(P)=1$ then $P\in\ass(M)$ if
and only if $P\supset J([2f])=\langle 2f\rangle$,
so $P=\langle 2\rangle$ or $\langle f\rangle$.
Notice that $J(\phi_2)=\langle 1\rangle$ so there are
no further primes, so $\ass(M)=\{\langle 2\rangle,\langle f\rangle\}$.

\vskip.1in
\noindent(c) The simplest setting in which a higher Fitting
ideal appears is Ledrappier's example.
Let $M=R/\langle 2,1+u_1+u_2\rangle$. A simple
syzygy calculation gives the free resolution
$$
0\longrightarrow R\overset{\phi_2}
{\longrightarrow}R^2\overset{\phi_1}
{\longrightarrow}R\longrightarrow M\longrightarrow0,
$$
where $\phi_1=\bmatrix 1+u_1+u_2,2\endbmatrix$ and $\phi_2=
\bmatrix 1+u_1+u_2\\-2\endbmatrix.$
If $\dim(P)=1$ then $P\in\ass(M)$ if
and only if $P\supset J([1+u_1+u_2,2])=\langle 2,1+u_1+u_2\rangle$,
so there are no primes here.
If $\dim(P)=2$ then $P\in\ass(M)$ if
and only if $P\supset J(\bmatrix 1+u_1+u_2\\-2\endbmatrix)
=\langle 2,1+u_1+u_2\rangle$, giving the one associated prime
$\langle 2,1+u_1+u_2\rangle$.

\vskip.1in
\noindent(d) 
Let $A=\bmatrix2&u_2^2-5&0\\0&u_1u_2-7u_1+u_2&3\endbmatrix$.
Then the first
Fitting ideal $J(A)$ is generated by the set $\{2u_1u_2-14u_1
+2u_2,6,3u_2^2-15\}$. A
principal prime ideal which contains $J(A)$ must contain $6$,
and must therefore be
generated by $2$ or $3$: in either case it cannot contain the other two
generators of $J(A)$. This proves that no principal ideals are associated to
the module $R^2/AR^3$. Using the special form of
the matrix we see that
the kernel of $A$ in $R^3$ is generated by the vector 
$$v=\bmatrix 3u_2^2-15\\-6\\2u_1u_2-14u_1+2u_2\endbmatrix.$$
The second Fitting
ideal $J(v)$ is equal to the first. Assume $P$ is prime with $\dim(R^P)=2$
and $P\supset J(v)$. Then this prime contains either $2$ or $3$.
If $3\in P$ then $P$ lies above the prime $\langle 3,u_1u_2-7u_1+2u_2\rangle$
which is the only one with $\dim(R^P)=2$. For the case $2\in P$ we have also
$u_2^2-5\in P$ but this element is modulo $2$ congruent to $(u_2-1)^2$;
this
means that
the only prime with the correct local dimension containing $2$
is $P=\langle 2,u_2-1\rangle$. The only associated primes
of $M=R^2/AR^3$ are
therefore $P_1=\langle 3,u_1u_2-7u_1+2u_2\rangle$ and $P_2=\langle 2,u_2-1
\rangle$. The corresponding dynamical system is expansive
and ergodic but not
mixing, and has zero entropy.

\vskip.1in
\noindent(e) 
Let $A=\bmatrix2&3u_2+5&3u_1-3u_2\\u_1-4&u_1-1&3u_1-6\endbmatrix$. Then the 
first Fitting ideal is generated by 
\begin{gather*}
  -3u_1+18-3u_1u_2+12u_2, \\
  18u_1-12-3u_1^2+3u_1u_2-12u_2\mbox{ and}\\
  -21u_2-30-3u_1^2+18u_1+12u_1u_2.
\end{gather*}
The only principal ideal
above $J(A)$ is $\langle3\rangle$. With a computer algebra system one can
calculate the kernel of the map $A$: it is generated by the vector 
$$
v=\bmatrix -7u_2-10-u_1^2+6u_1+4u_1u_2\\4+u_1^2-6u_1-u_1u_2+4u_2\\6-u_1u_2
+4u_2-u_1\endbmatrix.
$$
The second Fitting ideal $J(v)$ is generated by the
components of this vector and one
can calculate that
$$\{u_1-3u_2-4,3u_2^2+3u_2-2\}$$
is also a generating
set. So $J(v)$ is a prime with local dimension $2$. The only associated
primes of the module $R^2/AR^3$ are $\langle3\rangle$ and $\langle
u_1-3u_2-4,3u_2^2+3u_2-2\rangle$. The entropy of the corresponding dynamical
system is $\log3$, but the system does not have completely positive entropy. 
The ring $R/J(v)$ is isomorphic to a subring of ${\Bbb Q}
[\sqrt{ \frac{11}{12}}]$ via the map sending $u_2$ to
$-\frac{1}{2}+\sqrt{\frac{11}{12}}$ (a root of
$3y^2+y-2=0$), and $u_1$ to $\frac{5}{2}+\sqrt{\frac{33}{4}}$.
The field-theoretic norms of those two elements are $-2$
and $-\frac{2}{3}$ respectively. It follows that there can be no nontrivial 
$(n_1,n_2)\in{\Bbb Z}^2$ such that $u_1^{n_1}u_2^{n_2}-1\in J(v)$ because
this would yield $2^{n_1}(\frac{2}{3})^{n_2}=1$. The dynamical
system is therefore mixing of all orders and ergodic.

\vskip.1in
\noindent(f) Even in the square ($n=k$) case the first Fitting
ideal does not contain enough information to construct a
prime filtration of the module. The following type of
example is well--known (see for example
Remark 6(5) in \cite{ward-bernoulli-1992} or 
Example 5.3(2) in \cite{schmidt-1995-algebraic-dynamical-systems}).
Let $A=\bmatrix4-u_1&1\\1&-u_1\endbmatrix$ and
$B=\bmatrix3-u_1&2\\2&1-u_1\endbmatrix$. Then
$\det(A)=\det(B)$ so the systems $\alpha^{M_A}$
and $\alpha^{M_B}$ have the same entropy, number
of periodic points, and in fact are both isomorphic to
Bernoulli shifts and hence measurably isomorphic.
Both modules $M_A$ and $M_B$ have similar
finite free resolutions,
$$
0\longrightarrow R^2\overset{\phi}{\longrightarrow}
R^2\longrightarrow M\longrightarrow0,
$$
where $\phi=A$ for $M=M_A$ and
$\phi=B$ for $M=M_B$.
It is easy to check that $M_A\cong R/\langle
u_1^2-4u_1-1\rangle$, so that
$$
M_A\supset0
$$
is a prime filtration.
On the other hand, the shortest filtration of
$M_B$ is of the form
$$
M_B\supset  N\supset0,
$$
with first quotient $  N/\{0\}\cong R/\langle
u_1^2-4u_1-1\rangle$, and second quotient
$$M/  N\cong R/\langle
u_1^2-4u_1-1\rangle+Q$$
for some
ideal $Q\not\subset\langle
u_1^2-4u_1-1\rangle.$

\vskip.1in
\noindent(g) Let $f,g,h$ be co--prime elements of $R$,
and consider the module $M=R/[f,g,h]R^3$.
Then a finite free resolution is given by
the Koszul complex
$$
0\longrightarrow R\overset{\phi_1}
{\longrightarrow}R^3\overset{\phi_2}
{\longrightarrow}R^3\overset{\phi_3}\longrightarrow R\longrightarrow
M\longrightarrow0,
$$
in which $\phi_1=\bmatrix f\\g\\h\endbmatrix$,
$\phi_2=\bmatrix 0&h&-g\\-h&0&f\\g&-f&0\endbmatrix$
and $\phi_3=[f,g,h]$.

\vskip.1in
\noindent(h) An example in which the rank of the presenting
matrix is too small is given by $A=\bmatrix2\\1+u_1+u_2\endbmatrix.$
Let $M=M_A$; then
$\ann_M\left(\matrix1\\0\endmatrix\right)=
\langle0\rangle$, so $M$ has a free
submodule $L=\left(\matrix1\\0\endmatrix\right)R$.
The corresponding dynamical system therefore has as a
factor the full shift with circle alphabet, so
$h(\alpha^M)=\infty$.
The quotient $M/L\cong R/\langle1+u_1+u_2\rangle$ is then
of the form (\ref{finitepresentation}).
Of course the free submodule sits inside $M$ in many different
ways, so there is no ``canonical'' quotient $M/L$.

\vskip.1in
\noindent(i) The simplest examples
of algebraic dynamical systems without finite
presentation are certain non--expansive automorphisms of
solenoids, as studied in
\cite{lind-ward-1988} and
\cite{chothi-everest-ward-1997}.
Let $X=\widehat{\Bbb Z[\frac{1}{6}]}$, and let $\alpha$ be
the automorphism of $X$ dual to $x\mapsto 2x$
on $\Bbb Z[\frac{1}{6}]$
(here $d=1$). The $R$--module 
corresponding to the dynamical system
then has a chain of submodules
$$\textstyle{
\Bbb Z[\frac{1}{2}]\subset\frac{1}{3}\Bbb Z[\frac{1}{2}]\subset
\frac{1}{9}\Bbb Z[\frac{1}{2}]\subset\dots,}
$$
each of which is isomorphic as an $R$--module
to $R/\langle u_1-2\rangle$, that never stabilizes.
It follows that the corresponding module is not Noetherian.
\end{example}


\begin{thebibliography}{10}

\bibitem{chothi-everest-ward-1997}
V.~Chothi, G.~Everest and T.~Ward.
\newblock S-integer dynamical systems: periodic points.
\newblock {\em Journal f{\"u}r die Riene und. Angew.}, 489:99-132, 1997.

\bibitem{einsiedler-1997}
M.~Einsiedler.
\newblock A generalisation of {M}ahler measure and its application in algebraic
  dynamical systems.
\newblock {\em Acta Arithmetica, to appear}.

\bibitem{eisenbud}
D.~Eisenbud.
\newblock {\em Commutative Algebra with a view toward Algebraic Geometry}.
\newblock Springer Verlag, New York, 1995.

\bibitem{kitchens-schmidt-1989-main}
B.~Kitchens and K.~Schmidt.
\newblock Automorphisms of compact groups.
\newblock {\em Ergodic Theory and Dynamical Systems}, 9:691--735, 1989.

\bibitem{lind-ward-1988}
D.A~Lind and T.~Ward.
\newblock Automorphisms of solenoids and p-adic entropy.
\newblock {\em Ergodic Theory and Dynamical Systems}, 8:411--419, 1988.

\bibitem{lind-schmidt-ward-1990}
D.A. Lind, K.~Schmidt, and T.~Ward.
\newblock {M}ahler measure and entropy for commuting automorphisms of compact
  groups.
\newblock {\em Inventiones Math.}, 101:593--629, 1990.

\bibitem{matsumura-1986-commutative}
H.~Matsumura.
\newblock {\em Commutative Ring Theory}.
\newblock Cambridge University Press, Cambridge, 1986.

\bibitem{rudolph-schmidt-1995}
D.J. Rudolph and K.~Schmidt.
\newblock Almost block independence and {B}ernoullicity of {$Z^d$} actions by
  automorphisms of compact abelian groups.
\newblock {\em Inventiones Math.}, 120:455--488, 1995.

\bibitem{schmidt-1989}
K.~Schmidt.
\newblock Mixing automorphisms of compact groups and a theorem by {K}urt {M}ahler.
\newblock {\em Pacific Journal of Math.}, 137:371-385, 1989.

\bibitem{schmidt-1990}
K.~Schmidt.
\newblock Automorphisms of compact abelian groups and affine varieties.
\newblock {\em Proceedings of the London Math. Soc.}, 61:480--496, 1990.

\bibitem{schmidt-1995-algebraic-dynamical-systems}
K.~Schmidt.
\newblock {\em Dynamical Systems of Algebraic Origin}.
\newblock Birkh{\"a}user, Basel, 1995.

\bibitem{schmidt-ward-1993-schlickewei}
K.~Schmidt and T.~Ward.
\newblock Mixing automorphisms of compact groups and a theorem of
  {S}chlickewei.
\newblock {\em Inventiones Math.}, 111:69--76, 1993.

\bibitem{ward-bernoulli-1992}
T.~Ward.
\newblock The {B}ernoulli property for expansive {$Z^2$} actions on compact
  groups.
\newblock {\em Israel Journal of Math.}, 79:225--249, 1992.

\end{thebibliography}
\end{document}